\documentclass[12pt]{amsart}
\usepackage{epsfig,amssymb, eucal,xypic, url, amsthm}
\usepackage{color}

\newtheorem{theorem}{Theorem}[section]
\newtheorem{lemma}{Lemma}[section]

\newtheorem{proposition}{Proposition}[section]
\newtheorem{definition}{Definition}[section]
\theoremstyle{definition} \newtheorem{remark}{Remark}[section]
\newtheorem{example}{Example}[section]

\newcommand{\D}{{\mathcal D}}
\renewcommand{\S}{{\mathcal S}}

\renewcommand{\part}{\vdash}

\newcommand{\Hom}{\operatorname{Hom}}

\renewcommand{\S}{{\mathcal S}}

\renewcommand{\part}{\vdash}

\begin{document}
\title{Braid groups and Kleinian singularities}
\author{Christopher Brav and Hugh Thomas}

\begin{abstract}
We establish faithfulness of braid group actions generated by twists along an ADE configuration of $2$-spherical objects in a derived category. Our major tool is the Garside structure on braid groups of type ADE. This faithfulness result  provides the missing ingredient in Bridgeland's description of a space of stability conditions associated to a Kleinian singularity. 
\end{abstract}

\maketitle

\section{{\bf Introduction}}

The homological mirror symmetry program of Kontsevich \cite{homomirror} proposes a duality between symplectic geometry and complex geometry in the form 
of an equivalence between the derived Fukaya category on one side and the derived category of coherent sheaves on the other. Seidel and Thomas \cite{st} observed that since generalized Dehn twists around Lagrangian spheres in a symplectic manifold induce autoequivalences of the 
derived Fukaya category, there should be corresponding autoequivalences of the derived category of coherent sheaves on a mirror dual variety. 
With this expectation, they developed a theory of such autoequivalences, which they named spherical twists. 

In particular, they considered the case of a Kleinian singularity, constructed for example
as the quotient $\mathbb{C}^{2}/G$, where $G \subset SL_{2}(\mathbb{C})$ is a finite group. Any smoothing of the singularity is a symplectic manifold containing a 
collection of Lagrangian spheres in a configuration whose dual graph is a Dynkin diagram of type ADE. The generalized Dehn twists along these spheres are known to 
satisfy braid relations of type ADE in the symplectic mapping class group. 

A general pattern suggests that the minimal resolution $\pi: X \rightarrow \mathbb{C}^{2}/G$ 
should be mirror dual to the smoothing of $\mathbb{C}^{2}/G$, and so we should expect
to find spherical objects in $D^{b}(X)$ whose associated twists generate an action of a braid group on $D^{b}(X)$ by autoequivalences. Since the exceptional divisor $E=\pi^{-1}(0)$ consists of a tree of $-2$-curves whose dual graph is of type ADE, we expect that the desired spherical objects should come from the exceptional divisor. Indeed, the structure sheaves of the $-2$-curves are easily seen to be spherical and the associated twists are seen to satisfy braid relations in the group of autoequivalences up to isomorphism. In type A, using results from Khovanov and Seidel \cite{ks}, Seidel and Thomas were able to show that this action of the braid group is faithful.

Later, Thomas \cite{thomas} and Ishii-Ueda-Uehara \cite{iuu} (type A) and Bridgeland \cite{klein} (types ADE) studied the spaces of stability conditions of certain triangulated subcategories $\D$ of $D^{b}(X)$. Bridgeland showed that a connected component ${\rm Stab}_{0}(\D)$ of the space of stability conditions is a covering space of $\frak{h}^{reg}/W$, the space of regular orbits of the Weyl group $W$  corresponding to the singularity type, and that the braid group action on the derived category induces the full group of deck transformations of the cover. Moreover, he showed that faithfulness of the braid group action on the derived category implies faithfulness on 
${\rm Stab}_{0}(\D)$, and since $\pi_{1}(\frak{h}^{reg}/W)$ is known to be the braid group, such faithfulness implies that ${\rm Stab}_{0}(\D)$ is simply connected. Given 
the faithfulness result of \cite{st} in type A, simply-connectedness for ${\rm Stab}_{0}(\D)$ in type A follows.  

Our main goal is to provide the necessary faithfulness result to complete Bridgeland's description of spaces of stability conditions associated to Kleinian singularities in all types.
Specifically, we prove faithfulness for braid group actions in types ADE generated by twists along $2$-spherical objects (Theorem \ref{faithful}). 

Our proof makes essential use of the Garside structure on braid groups. We expect that similar methods can be used to study other braid group actions on categories appearing in algebraic geometry and representation theory. 

We summarize the contents of the paper.
In Section 2, we review some of the theory of braid groups to establish notation. In Section 3, we review the theory of  spherical twists, and then prove our main result, Theorem \ref{faithful}: faithfulness of braid group actions generated by $2$-spherical twists. In Section 4, we recall Bridgeland's work concerning spaces of stability conditions associated to Kleinian singularities,
pointing out how Theorem \ref{faithful} applies to this situation.
\\

\noindent{\bf Acknowledgements}
We are grateful to Tom Bridgeland, Sabin Cautis, Colin Ingalls, Joel Kamnitzer, Bernhard Keller, David Ploog, Leonid Positselski, and Mike Roth for helpful conversations about this project. We wish to express particular gratitude to an anonymous referee who pointed out an error in our attempted proof of faithfulness of extended affine braid group actions, which had been included in an earlier version of this paper. HT was partially supported by an NSERC discovery grant.

\section{{\bf Background on braid groups}}
We begin by establishing notation for Weyl groups, braid groups and related structures associated to a Dynkin diagram $\Gamma$ and summarizing what we shall later need to know about Garside factorizations of  braid group elements. A good general reference for this material is Kassel-Turaev, \cite[Chapter 6]{kassel}.

The Weyl group and braid group associated to a Dynkin diagram $\Gamma$ have various
geometric, topological, and combinatorial realizations. For our purposes, however, it will be sufficient to describe these groups by presentations.

\medskip
\noindent{\bf Weyl groups and braid groups}
\medskip

\begin{definition} Given a Dynkin diagram $\Gamma$ of type {\rm ADE},
the associated Weyl group $W$ has generators $s_{i}$ with $i \in \Gamma$ a node, subject to the relations $s_{i}^{2}=1$, $s_{i}s_{j}=s_{j}s_{i}$ if $i,j$ are not adjacent in $\Gamma$, and $s_{i}s_{j}s_{i}=s_{j}s_{i}s_{j}$ if $i,j$ are not adjacent in $\Gamma$.
\end{definition}

Given $w \in W$, we call an expression $w=s_{i_{1}}\cdots s_{i_{k}}$ {\it reduced} if there are no shorter expressions of $w$ in terms of the generators. Now define a length function $\ell: W \rightarrow \mathbb{N}$, where
$\ell(w)$ is the length of a reduced expression. We say that a factorization $w=uv$ is reduced if 
$\ell(w)=\ell(u)+\ell(v)$. Given a reduced factorization $w=uv$, we say that $u$ is a left factor of $w$ and $v$ a right factor.
It is known that there is a unique longest element $w_{0} \in W$ and that every element
$w \in W$ is a left factor and a right factor of $w_{0}$.

\begin{definition}
The braid group $B$ is generated by $\sigma_{i}, i \in \Gamma$, 
subject to the braid relations 
$\sigma_{i}\sigma_{j}=\sigma_{j}\sigma_{i}$ if $i,j$ are not adjacent in $\Gamma$ and 
$\sigma_{i}\sigma_{j}\sigma_{i}=\sigma_{j}\sigma_{i}\sigma_{j}$ if $i,j$ are adjacent in $\Gamma$.

The braid monoid $B^{+}$ is given by the same presentation, but now in the category of monoids.
\end{definition}

It is known that the natural monoid homomorphism $B^{+} \rightarrow B$
sending $\sigma_{i}$ to $\sigma_{i}$ is an injection and identifies $B$ with the group of fractions of $B^{+}$. That is to say, any monoid homomorphism $\rho^{+}: B^{+} \rightarrow G$ to a group $G$ extends uniquely to a group homomorphism $\rho: B \rightarrow G$. 
Note that the image of $B^{+} \rightarrow B$ is just the submonoid of $B$ generated by the $\sigma_{i}$.

\medskip

\begin{remark}\label{fundamental}
It is well-known that the braid group $B$ of type $\Gamma$
may be realized as the fundamental group
of the space of regular orbits $\frak{h}^{reg}/W$ for the Weyl group $W$ acting on the the Cartan algebra $\frak{h}$ of type $\Gamma$.  
\end{remark}

Since the relations among the generators $\sigma_{i}$ of $B$ are also satisfied by the generators $s_{i}$ of the associated Weyl group $W$, we have a natural surjection $\pi: B \rightarrow W$ sending $\sigma_{i}$ to $s_{i}$. There is, moreover, a set-theoretic
section $\varphi: W \rightarrow B$ of $\pi$, which sends $s_{i}$ to $\sigma_{i}$ and sends
an element $w \in W$ to the product of $\sigma_{i}$ corresponding
to a reduced expression in terms of the $s_{i}$. This prescription is well-defined since any two reduced expressions of an element $w \in W$ can be related by braid relations, which also hold in $B$. In particular, let $\Delta \in B$ be the image of the longest word $w_{0}$ 
under $\varphi$. Note that the image of the section $\varphi(W)$ is by construction
contained in the braid monoid $B^{+}$. For brevity,
we shall often denote $\varphi(w)$ by $\tilde{w}$. Note that $\widetilde{uv}=\tilde{u}\tilde{v}$ if $\ell(uv)=\ell(u)+\ell(v)$.

On the braid monoid $B^{+}$, we have a length function $\ell: B^{+} 
\rightarrow \mathbb{N}$, again defined as the length of a shortest expression of an element in terms of the generators. Note that $\alpha \in \varphi(W)$ if and only 
if $\ell(\alpha)=\ell(\pi(\alpha))$. In particular,
the length of $w \in W$ is the same as the length of $\widetilde{w} \in B^{+}$.

A factorization $\alpha=\beta\gamma$ in $ B^+$ is said to be reduced 
if $\ell(\alpha)=\ell(\beta) + \ell(\gamma)$. Given a reduced 
factorization
$\alpha=\beta\gamma$ in $ B^+$,  we say that $\beta$ is a left factor of $\alpha$ and $\gamma$ a right factor. Note that the image $\varphi(W)$ can be described as the 
set of left factors or the set of right factors of $\Delta$.  In particular, $\varphi(W)$ is
closed under taking left or right factors. 

\medskip
\noindent{\bf Garside factorization}

\medskip
We now want to describe a normal form, the {\it Garside factorization}, for elements of the braid group $B$. For a more thorough discussion and further references, see \cite[Chapter 6]{kassel}.  

\medskip 

The Garside factorization for an element $\alpha \in B^+$ is 
a reduced expression $\alpha=\alpha_{k}\cdots \alpha_{1}$, where $\alpha_{i}=\widetilde{w}_{i} \in \varphi(W)$ are canonically defined elements of the braid monoid. Let $\alpha \in B^+$.    
It can be shown that $\alpha$ has a unique longest right factor 
lying in $\varphi(W)$ and that all other right factors are right factors of the longest one. By definition, we take this longest right factor to be $\alpha_{1}$. Writing $\alpha=\alpha'\alpha_1$, the 
succeeding factors are defined by applying the same procedure recursively 
to $\alpha'$.  

The following standard result says that the property of being a Garside
factorization can be checked locally.  

\begin{lemma} For $\alpha_i\in B^+$, $1\leq i \leq k$, 
$(\alpha_k,\dots,\alpha_1)$ is the Garside factorization for
$\alpha_k\dots\alpha_1$ if and only if the Garside factorization for $\alpha_i\alpha_{i-1}$
is $(\alpha_i,\alpha_{i-1})$.  \end{lemma}

The previous lemma combines well with the following, which gives a
more explicit procedure for checking that $(\alpha,\beta)$ is a
Garside factorization of $\alpha\beta$.

\begin{lemma}\label{checkgarside}
$(\tilde{u},\tilde{v})$ is the Garside factorization of 
$\tilde{u}\tilde{v}$ if and only if for any $s_i$ which can appear as the rightmost factor
of $u$, $s_i$ can also appear as the leftmost factor of $v$.  
\end{lemma}

\begin{proof}
Suppose first that $(\tilde{u},\tilde{v})$ is the Garside factorization of
$\tilde{u}\tilde{v}$.
Consider any $s_i$ which can appear as a rightmost factor of $u$.  If
$\ell(s_iv)=\ell(v)+1$, then $(\widetilde{us_i},\widetilde{s_iv})$ would be a 
factorization of $\tilde{u}\tilde{v}$ with longer right factor, contradicting
our assumption.  Thus $\ell(s_iv)<\ell(v)$, which implies that $v$ can be written with $s_i$ as its leftmost factor.  

Conversely, assume that any $s_{i}$ appearing as a right factor of $u$ 
also appears as a left factor of $v$. If $\tilde{u}\tilde{v}$ is 
not already the Garside factorization, then let $\tilde{v}'$ be the 
rightmost Garside factor. Then $v'$ must have $v$ as a 
right factor and there is some reduced expression $v'=ws_iv$, so 
$\ell(s_iv)>\ell(v)$.    
Then we can write $\tilde{u}\tilde{v}=
\alpha\sigma_i \tilde{v}$. Thus $\alpha\sigma_i=\widetilde u$.  It follows that $\pi(\alpha)s_i$ is a reduced factorization of $u$, so $u$ has a factorization $u=xs_i$. Since $u$ has $s_i$ as a right factor, but
$\ell(s_iv)>\ell(v)$, $v$ does not have $s_{i}$ as a left factor, 
contrary to assumption.
\end{proof}

More generally, for $\alpha \in B$, the Garside factorization is of the form
$(\alpha_k,\dots,\alpha_1)$ where $\alpha_i\in \varphi(W)\cup \{\Delta^{-1}\}$.
More precisely, there is some $j$ such that:
\begin{itemize}
\item 
$\alpha_1=\dots=\alpha_j=\Delta$ or $\alpha_1=\dots=\alpha_j=\Delta^{-1}$.
\item
For $i>j$, $\alpha_i\in \varphi(W)\setminus\{\Delta\}$.
\item $(\alpha_k,\dots,\alpha_{j+1})$ is the Garside factorization of
$\alpha_k\dots\alpha_{j+1}$ in $B^+$.  
\end{itemize}
Further, any factorization satisfying these properties is the Garside
factorization of some element of $B$.  

\medskip

The following lemma will be useful in establishing faithfulness of braid group actions.

\begin{lemma}\label{monoid}
A group homomorphism $\rho: B \rightarrow G$ is injective if and only if the induced monoid homomorphism $\rho^{+}: B^{+} \hookrightarrow B \rightarrow G$
is injective.
\end{lemma}

\begin{proof}
Injectivity of $\rho$ clearly implies injectivity of $\rho^{+}$. Conversely, suppose $\rho^{+}$ is injective and let $\rho(\alpha)=1$ for $\alpha \in B$.
Using the Garside factorization, write $\alpha=\beta\Delta^{m}$ for $\beta \in B^{+}$ and $m \in \mathbb{Z}$. If $m \geq 0$, then $\alpha \in B^{+}$
and so $\rho(\alpha)=\rho^{+}(\alpha)=1$ implies $\alpha=1$. If $m < 0$, then $\rho(\alpha)=1$ gives $\rho^{+}(\beta)=\rho^{+}(\Delta^{-m})$, so injectivity of $\rho^{+}$ 
implies $\beta=\Delta^{-m}$ and hence $\alpha=1$.
\end{proof}

\section{{\bf Spherical twists and braid group actions}}

We begin by reviewing some aspects of the theory of spherical objects (introduced in \cite{st} and refined in \cite{raph} and \cite{rina}). Since it is sufficient for our 
applications and the statement of our main result, we restrict to the two dimensional case to simplify our exposition.

\medskip

Throughout this paper, any triangulated category
$\D$ is assumed to be linear over a field $\mathrm{k}$ and to come with a fixed enhancement. For instance, we may take $\D$ to be the homotopy category
of a stable $\infty$-category in the sense of \cite{lurie}, an algebraic triangulated category in the sense of \cite{kellertilting}, or a pre-triangulated differential-graded category in the sense of \cite{bonkap}.

Since $\D$ comes with a fixed enhancement, we have functorial cones, 
derived Hom-complexes ${\rm RHom}(X,Y)$ for any two objects $X,Y \in 
\D$, and the adjoint pair of functors $? \otimes X: D(\mathrm{k}) 
\rightarrow \D$ and ${\rm RHom}(X,?): \D \rightarrow D(\mathrm{k})$. 
If $X$ is such that ${\rm RHom}(M,X)$ has total finite dimensional cohomology for all $M \in \D$, then the functor $?\otimes X : D(\mathrm{k}) \rightarrow \D$ also has a left adjoint ${\rm RHom}(?,X)^{\vee}: \D \rightarrow D^{b}(\mathrm{k})$, where
$\vee$ denotes the dualization ${\rm RHom}(?, \mathrm{k}): D^{b}(\mathrm{k}) \rightarrow D^{b}(\mathrm{k})$. (Here $D(\mathrm{k})$ is the unbounded derived category of vector spaces and $D^{b}(\mathrm{k})$ is the bounded derived category of finite dimensional $\mathrm{k}$-vector spaces.)

In order to simplify notation, we often write $[X,Y]$ for ${\rm Hom}_{\D}(X,Y)$, $[X,Y]_{d}$ or ${\rm Ext}^{d}(X,Y)$ for ${\rm Hom}(X,Y[d])$, and $[X,Y]_{*}$ for
$\bigoplus_{d} {\rm Hom}(X,Y[d])$.  

For brevity, let us temporarily denote
${\rm RHom}(X,?)$ by $R$, $?\otimes X$ by $X$, and ${\rm RHom}(?,X)^{\vee}$ by $L$.
Then we have the following units and counits:
\begin{equation}\label{unit}
{\rm Id}_{\mathrm{k}}   \rightarrow RX
\end{equation}
\begin{equation}\label{counit}
XR \rightarrow {\rm Id}_{\D}\end{equation}
\begin{equation}\label{unit2}
{\rm Id}_{\D} \rightarrow XL
\end{equation}
\begin{equation}\label{counit2}
LX \rightarrow
{\rm Id}_{\mathrm{k}}
\end{equation}

Now define an exact endofunctor $\Phi_{X}$ of 
$D(\mathrm{k})$ as the cone of the unit in (\ref{unit}), so that we have a triangle of functors ${\rm Id}_{\mathrm{k}}   \rightarrow RX \rightarrow \Phi_{X}
$. Applying $R$ on the left of (\ref{unit2}), we get a morphism $R \rightarrow RXL$, and applying $L$ to the right of $RX \rightarrow \Phi_{X}$, 
we get a morphism $RXL \rightarrow \Phi_{X} L$. Composing these two morphisms gives a morphism
\begin{equation} \label{duality}
R  \rightarrow \Phi_{X} L.
\end{equation}

\medskip

\begin{definition}
Let $S \in \D$ and suppose that
${\rm RHom}(?,S)^{\vee}$ always has total finite dimensional cohomology so that it gives a left adjoint to the functor $?\otimes S$. Then $S$ is {\rm 
$2$-spherical} if 
\begin{enumerate}
\item the cone $\Phi_{S}$ of ${\rm Id}_{\mathrm{k}} \rightarrow {\rm RHom}(S,?\otimes S)$ is isomorphic to the shift functor $[2]$
\item the natural morphism ${\rm RHom}(S,?) \rightarrow {\rm RHom}(?,S[2])^{\vee}$ from
(\ref{duality}) is an isomorphism.
\end{enumerate}
\end{definition}

\begin{remark} Applying the first condition in the above 
definition to the object $\mathrm{k}$, we see that ${\rm Ext}^{*}(S,S)$ 
is one-dimensional in degrees $0$ and $2$ and zero elsewhere, hence 
isomorphic to the cohomology of the $2$-sphere 
(hence the name `spherical object'). The second condition is a 
Calabi-Yau condition and says that 
the shift functor $[2]$ restricted to $S$ realizes a kind of Serre 
duality. 

Note that condition b) implies that composition of morphisms gives a perfect pairing
\begin{equation}\label{pairing}
[S, X] \otimes [X,S]_{2} \rightarrow [S,S]_{2}\simeq \mathrm{k}.
\end{equation}
\end{remark}

\medskip

\begin{example} Let $X$ be a smooth quasi-projective surface, $C \subset X$ a $-2$-curve. Then any twisted structure sheaf $\mathcal{O}_{C}(d)$ is 
$2$-spherical in $D^{b}(X)$, the bounded derived category of coherent sheaves on $X$.
\end{example}

\medskip

We shall be particularly interested in special configurations of spherical objects.

\begin{definition}
Let $\Gamma$ be a Dynkin diagram of type ADE. We say that
a collection of $2$-spherical objects $S_{i}$, $i \in \Gamma$, is a 
$\Gamma$-configuration if for $i \neq j$, the space $[S_{i},S_{j}]_{*}$ 
is one-dimensional and concentrated in degree $1$  when $i$ and $j$ are 
adjacent in $\Gamma$ and is zero otherwise.
\end{definition}

\begin{example}\label{kleinian}
Consider a finite subgroup $G \subset SL_{2}(\mathbb{C})$ and the quotient $\mathbb{C}^{2}/G$, a singular surface with a unique singular point $0$, commonly called a Kleinian singularity after Felix Klein who first determined the ring of invariants $\mathbb{C}[x,y]^{G}$ in his classic book The Icosahedron \cite{kleinf}. The minimal resolution $\pi: X \rightarrow \mathbb{C}^{2}/G$ was later studied by du Val \cite{duval}, who showed that the exceptional divisor $E=\pi^{-1}(0)$ consists
of a tree of $-2$-curves whose dual graph $\Gamma$ is a Coxeter-Dynkin diagram $\Gamma$ of type ADE. 
The geometry of the minimal resolution $\pi: X \rightarrow \mathbb{C}^{2}/G$
thus provides a beautiful bijection between (conjugacy classes of) 
finite subgroups $G \subset SL_{2}(\mathbb{C})$ and Dynkin diagrams of 
type ADE.

Letting $E_{i}$, $i \in \Gamma$ be the irreducible components of $\pi^{-1}(0)=E$, a standard computation shows that  
the collection $S_i=\mathcal{O}_{E_i}$ is a $\Gamma$-configuration of $2$-spherical objects. (The same is of course true for the collection $S_i=\mathcal{O}_{E_i}(-1)[1]$, which, despite its appearance, is sometimes more convenient to work with.)
\end{example}

\medskip

Returning to a general $\Gamma$-configuration, we see that 
since the objects $S_{i}$ are $2$-spherical, we have perfect pairings as in (\ref{pairing}):  
$$[S_{i},S_{j}]_{1}\otimes [S_{j},S_{i}]_{1} \rightarrow [S_{i},S_{i}]_{2}\simeq \mathrm{k}.$$
As a consequence, we have the following lemma that will be useful in the proof of Proposition \ref{probe} below.
 
\begin{lemma}\label{factor}
If $i$ and $j$ are adjacent in $\Gamma$, then any morphism $S_{i} \rightarrow S_{i}[2]$ factors as 
$S_{i} \rightarrow S_{j}[1] \rightarrow S_{i}[2]$.
\end{lemma}

\medskip

We now recall some results of Seidel-Thomas \cite{st} concerning spherical 
twists in enhanced triangulated categories, stated for $2$-spherical objects for simplicity, and then give the main result 
of this paper, Theorem \ref{faithful}: braid group actions of types ADE, generated by
$2$-spherical twists , are faithful in types ADE.

\medskip

Given an object $S \in \D$, we have a functor
$? \otimes S: D(\mathrm{k}) \rightarrow \D$ and its right adjoint 
${\rm RHom}(S,?): \D \rightarrow D(\mathrm{k})$. We therefore can define a `twist' functor $t_{S}: \D \rightarrow \D$ as the cone
of the counit ${\rm RHom}(S,?) \otimes S \rightarrow {\rm Id}_{\D}$. When $S$ is a spherical object, Seidel and Thomas \cite{st}
showed that $t_{S}$ is an (auto)equivalence.

We collect some standard facts about spherical twists in the following lemma.
\begin{lemma}\label{twists}
\begin{enumerate}
\item If $S$ is $2$-spherical, $t_{S}(S) \simeq S[-1]$
\item If $F$ is an autoequivalence of $\D$, then $F \circ t_{S} \simeq t_{F(S)} \circ F$.
\item If $S_{i}$ and $S_{j}$ are not adjacent in a $\Gamma$-configuration of spherical objects, then $t_{i}S_{j}\simeq S_{j} $.
\item If $S_{i}$ and $S_{j}$ are adjacent in a $\Gamma$-configuration of spherical objects, then $t_{i}t_{j}S_{i} \simeq S_{j}$.
\end{enumerate}
\end{lemma}

From Lemma \ref{twists} c) and d), it can be shown that if $S_{i}, i \in \Gamma$ form a $\Gamma$-configuration, then the associated twist
functors $t_{i}:=t_{S_{i}}$ satisfy the braid relations of type $\Gamma$, up to isomomorphism, and so there is a homomorphism 
$$\rho: B \rightarrow {\rm Aut}(\D)$$ from the braid group $B$ of type $\Gamma$ to the group ${\rm Aut}(\D)$ of isomorphism classes of autoequivalences
of $\D$. We denote the functor $\rho(\alpha)$ by $t_{\alpha}$. 

In type $A$, Seidel and Thomas, based on work of Khovanov and Seidel \cite{ks}, showed that the homomorphism $\rho$ is injective. We will generalize this to types ADE.
Our stragey is to show that a braid group element $\alpha$ is completely determined by the action of the corresponding twist $t_{\alpha}$ on $\S=\bigoplus_{i} S_{i}$. 
To do this, we will probe $t_{\alpha}\S$ by considering the Hom spaces $[S_i, t_{\alpha}\S]_{k}$. The following proposition provides the necessary information.

\medskip

\begin{proposition}\label{probe}
Let $1\ne \alpha \in B^{+}$ have Garside factorization $\alpha=\tilde{w}_{k}\cdots \tilde{w}_{1}$. Then 
\begin{enumerate}
\item $[\S, t_{\alpha}\S]_{d}=0$ for $d > k+2$
\item $[S_{i}, t_{\alpha}\S]_{k+2}\neq 0$ if and only if $s_{i}$ is a left factor of $w_{k}$ (in particular $[\S, t_{\alpha}\S]_{k+2}\neq 0$ ).
\item The maximal $p$ such that $[\S,t_{\alpha}\S]_{p} \neq 0$ is precisely $p=k+2$.
\end{enumerate}
\end{proposition}

In words, we can determine the number $k$ of Garside factors of $\alpha$ from the maximal degree $p$ of a map from $\S$ to $t_{\alpha}\S$ and we can determine 
whether or not $s_{i}$ is a left factor of the final Garside factor of $\alpha$ by studying maps from $S_{i}$ to $t_{\alpha}\S$.

\medskip

Before proving Proposition \ref{probe}, let us see how it implies our main result.

\begin{theorem}\label{faithful}
The homomorphism $\rho: B \rightarrow {\rm Aut}(\D)$ is injective.
\end{theorem}

\begin{proof}
By Lemma \ref{monoid}, it is enough to show injectivity on $B^{+}$. To do this, we show that $\alpha \in B^{+}$ can be recovered
from the functor $t_{\alpha}$. Thus if two functors $t_{\alpha}$ and $t_{\beta}$ for $\alpha, \beta \in B^{+}$ are isomorphic,
then we must have $\alpha=\beta$.

To recover $\alpha$ from $t_{\alpha}$, we study the mapping space $[\S,t_{\alpha}\S]_{*}$. By Proposition \ref{probe} c),
we know that the number of Garside factors of $\alpha$ is $k=p-2$ where $p$ is the maximal degree of a non-zero map
from $\S$ to $t_{\alpha}\S$. Now let $\alpha=\tilde{w}_{k}\cdots \tilde{w}_{1}$ be the Garside factorization of $\alpha$
and let $\beta=\tilde{w}_{k-1}\cdots \tilde{w}_{1}$

First, we will determine
a reduced decomposition of $w_{k}$ and hence $w_{k}$ itself. Since $[\S, t_{\alpha}\S]_{k+2} \neq 0$, there must be some $S_{i}$
such that $[S_{i},t_{\alpha}\S]_{k+2} \neq 0$, and by Proposition \ref{probe} b), $s_{i}$ must then be a left factor of $w_{k}$, so 
write a reduced expression $w_{k}=s_{i}u$. Then $t_{i}^{-1}t_{\alpha}=t_{\tilde{u}\beta}$. Now consider $[\S, t_{\tilde{u}\beta} \S]_{k+2}$.
If it is zero, then $\tilde{u}\beta$ has $k-1$ Garside factors, so $u=1$ and we have determined that $w_{k}=s_{i}$. Otherwise
we repeat the above argument to find a left factor $s_{j}$ of $u$. Proceeding in this way, we eventually find a reduced
decomposition $w_{k}=s_{i}s_{j} \cdots$. 

Once we have determined $w_{k}$, we repeat the whole process on $t_{w_{k}}^{-1}t_{\alpha}=t_{\beta}$ to determine
$w_{k-1}$, and so on, until we have determined in order all of the Garside factors of $\alpha$ and hence $\alpha$ itself.
\end{proof}

We shall need the following lemma in the proof of Proposition \ref{probe}.

\begin{lemma}\label{max}
Let $Y \in \D$.
\begin{enumerate}

\item If $l$ is maximal such that $[S_{i},Y]_{l} \neq 0$, then $m=l+1$ is maximal such that  $[S_{i}, t_{i}Y]_{m}\neq 0$.

\item Let $p$ be maximal such that $[\S,Y]_{p} \neq 0$. If $q$ is maximal such that $[\S, t_{i}Y]_{q} \neq 0$, then $p \leq q \leq p+1$. Further, $q=p+1$
if and only if $[S_i,Y]_p \ne 0$.  
\end{enumerate}
\end{lemma}
In words, part a) says that twisting an object $Y$ by $t_{i}$ increases by one the maximal degree of a map from $S_{i}$.  Part b) says that twisting by $t_{i}$ cannot 
decrease the maximal degree of a map from $\S$, it increases the maximal degree if and only if there is a map of degree $p$ from
$S_i$ to $Y$, and, if so, it increases the maximal degree by one.  

\begin{proof}
\noindent a) Since $t_{i}S_{i}\simeq S_{i}[-1]$, twisting
a non-zero morphism $S_{i} \rightarrow Y[l]$ produces
a non-zero morphism $S_{i} \rightarrow t_{i}Y[l+1]$. It must be of maximal degree, since if there were a non-zero morphism $S_{i} \rightarrow t_{i}Y[m]$ with 
$m > l+1$, then twisting by $t_{i}^{-1}$ and translating would give a non-zero map $S_{i} \rightarrow Y[m-1]$
of degree greater than $l$, contrary to assumption.

\medskip

\noindent b) Consider the triangle
\begin{equation*}
\xymatrix{[S_{i},Y]_{*} \otimes S_{i} \ar[r]^{\;\;\;\;\;\;\;\;\; f}& Y \ar[r]^{g} & t_{i}Y}
\end{equation*}
Let $p$ be maximal such that $[\S,Y]_{p} \neq 0$ and $q$ be maximal such that $[\S,t_{i}Y]_{q} \neq 0$. For $j\in \Gamma$, let $l_j$ be the maximum degree of a map from $S_j$ to $Y$ and note that $l_{j} \leq p$ for all $j$. 
By part a), the maximum degree of a map from $S_i$ to $t_iY$ is $l_{i}+1$.

Now consider some $j \ne i$.
The maximum degree of a map from $S_j$ to 
$[S_i,Y]_*\otimes S_i[1]$ is $l_i$.  Since $t_iY$ lies in a triangle between
$Y$ and $[S_i,Y]_*\otimes S_i[1]$, the maximum degree of a map from
$S_j$ to $t_iY$ is no greater than $\max(l_i,l_j) \leq p+1$.

This shows that $q \leq p+1$, and that $q=p+1$ if and only if the maximum 
degree of a map from $S_i$ to $Y$ is $p$. 

To show that $q\geq p$, observe first that if $l_i\geq p-1$, we are done,
so assume that $l_i <p-1$.  Now, by assumption, for some $j\ne i$, 
there must be a nonzero map in $[S_j, Y]_p$.  Since $[S_i,Y]_{p-1}=0$,
we know by the condition of a $\Gamma$-configuration that $[S_j, [S_i,Y]_*\otimes S_i]_p=0$, so the nonzero map in
$[S_j,Y]_p$ must compose with $g$ to give a non-zero map in
$[S_j,t_iY]_p$.  This shows that $q\geq p$.  .
\end{proof}

We now prove Proposition \ref{probe}.

\begin{proof}[Proof of Proposition \ref{probe}]
We prove a) and b) together by simultaneous induction on $k$ and $\ell(w_k)$.
Statement c) follows immediately from a) and b).

The base case of the induction is when $k=1$ and $w_1=s_i$ for some $i \in \Gamma$.  
In this case, the statements follow from a straight-forward calculation or directly from Lemma \ref{max}.

Now we suppose that we know a) and b) for any $1\ne\beta \in B^+$ 
with fewer Garside factors
or with shorter final Garside factor than $\alpha$, 
and we prove a) and b) for $\alpha$.  

\medskip

\noindent a) Suppose first that $w_k=s_i$ for some $i$.  In this case,
we know by the induction hypothesis that 
$[\S,t_{\tilde w_{k-1}\dots \tilde w_1}\S]_d=0$ for $d > k+1$, so it follows
from Lemma \ref{max} that $[\S,t_\alpha \S]_d=0$ for $d>k+2$.

Now suppose that $\ell(w_k)>1$, and fix a reduced decomposition $w_k=s_iu$.  
Let $\beta=\tilde{u}\tilde{w}_{k-1}\dots\tilde{w}_1$.  By the induction
hypothesis, we know that $[\S,t_\beta\S]_d=0$ for $d>k+2$.  Since $s_iu$ is
reduced, we know that $u$ cannot be written with $s_i$ as a leftmost factor.
Therefore, by b), we know that $[S_i,t_\beta\S]_{k+2}=0$.  By Lemma \ref{max},
it follows that $[\S,t_it_\beta\S]_d=0$ for $d>k+2$, as desired.  

\medskip

\noindent b) First we prove that if $s_{i}$ is a left factor of $w_{k}$, then $[S_{i}, t_{\alpha}\S]_{k+2} \neq 0$. 
For brevity, write $\beta=\tilde{w}_{k-1} \cdots \tilde{w}_{1}$.

\medskip

Consider first the case that $\ell(w_k)=1$.  Since we have already disposed of the case that $k=1$ and $\ell(w_1)=1$, we may assume that $k>1$. By Lemma \ref{checkgarside}, $w_{k-1}$ must have $s_{i}$ as a left factor, so by induction on $k$, 
there is a non-zero map $S_{i} \rightarrow t_{\beta}\S[k+1]$. Again by Lemma \ref{max}, part a), twisting
with $t_{i}=t_{\tilde{w}_{k}}$ then produces a non-zero map
$S_{i} \rightarrow t_{\alpha}\S[k+2]$, as needed.   

\medskip

Now assume that $\ell(w_k)>1$ and suppose then that we have a reduced expression $w_{k}=s_{i}s_{j}u$. This leads us to 
consider various cases.

\medskip

\noindent {\bf Case 1} Suppose $s_{i}$ and $s_{j}$ commute (so $[S_{i},S_{j}]_{*}=0$) and thus we may write $w_{k}=s_{j}s_{i}u$. By induction 
on the length of $w_{k}$, we have a non-zero map $S_{i} \rightarrow t_{i}t_{\tilde{u}}t_{\beta}\S[k+2]$. Since $t_{j}S_{i} \simeq S_{i}$, twisting with 
$t_{j}$ produces a non-zero map $S_{i} \rightarrow t_{\alpha}\S[k+2]$, as desired. 

\medskip

\noindent {\bf Case 2} Suppose $s_{i}$ and $s_{j}$ do not commute and we 
have a reduced decomposition $w_{k}=s_{i}s_{j}s_{i}v=s_{j}s_{i}s_{j}v$. 
By induction on the length of $w_{k}$, we have a non-zero map 
$S_{j} \rightarrow t_{j}t_{\tilde v}t_{\beta}\S[k+2]$. 
Applying $t_{j}t_{i}$ gives a non-zero map $S_{i}\rightarrow 
t_{\alpha}\S[k+2]$ (since $t_{j}t_{i}S_{j} \simeq S_{i}$ by 
Lemma \ref{twists}). 

\medskip

\noindent {\bf Case 3} Suppose $s_{i}$ and $s_{j}$ do not commute and 
there is no reduced decomposition $w_{k}=s_{i}s_{j}s_{i}v$, so $s_{i}$ is {\it 
not} a left factor of $u$. (Note that this includes the case $u=1$).

We need to show that there is a non-zero map $S_{i}[-k-2] \rightarrow t_{i}t_{j}t_{\tilde{u}}t_{\beta}\S$. Applying $t_{i}^{-1}$, we see that this is 
equivalent to having a non-zero map
$S_{i} \rightarrow t_{j}t_{\tilde{u}}t_{\beta}\S[k+1]$ (since $t_{i}^{-1}S_{i}\simeq S_{i}[1]$). 

Now consider the triangle 
\begin{equation}\label{triangle}
[S_{j}, t_{\tilde{u}}t_{\beta}\S]_{*} \otimes S_{j}
 \rightarrow  t_{\tilde{u}}t_{\beta}\S \rightarrow 
t_{j}t_{\tilde{u}}t_{\beta}\S 
\end{equation}
Note $[S_{j}, t_{\tilde{u}}t_{\beta}\S]_{k+2} 
=0$, or else $s_{j}$ would be a left factor 
of $u$ (by induction on the length of $w_{k}$) and then the expression for 
$w_{k}$ would not be reduced, or $u=1$, and then we may apply 
induction on the number of factors. 

Now we  claim that $[S_{j}, t_{\tilde{u}}t_{\beta}\S]_{k+1} \neq 0$. Suppose otherwise.  By induction on length of $w_{k}$, there
is a non-zero map $S_{j}[-k-2] \rightarrow t_{j}t_{\tilde{u}}t_{\beta}\S$. 
But in the triangle (\ref{triangle}), the objects to the left and 
right of $t_{j}t_{\tilde{u}}t_{\beta}\S$ admit 
no map from $S_{j}[-k-2]$, a contradiction. Thus
$[S_{j}, t_{\tilde{u}}t_{\beta}\S]_{k+1} \neq 0$.

Then since $i$ and $j$ are neighbors, we get a map
$\varphi: S_{i}[-k-2] \rightarrow [S_{j}, t_{\tilde{u}}t_{\beta}\S]_{k+1} \otimes S_{j}[-k-1]$. Considering again the triangle (\ref{triangle}),
we see that the composition to $t_{\tilde{u}}t_{\beta}\S$ vanishes, so $\varphi$ must factor through 
$t_{j}t_{\tilde{u}}t_{{\beta}}\S[-1]$, giving a non-zero
map $S_{i}[-k-2] \rightarrow t_{j}t_{\tilde{u}}t_{\beta}\S[-1]$. We twist this with $t_{i}$ and shift to get a 
non-zero map $S_{i}[-k-2] \rightarrow t_{\alpha}\S$, as desired.

\medskip

We now establish the opposite implication by proving the contrapositive: if $s_i$ is not a left factor of $w_k$, then 
$[S_i,t_\alpha \S]_{k+2}=0$.  

Fix a reduced factorization $w_k=s_ju$.  
Write $\beta=\tilde u\tilde{w}_{k-1} \cdots \tilde{w}_1$. 
Consider the triangle:

\begin{equation}\label{triangle2}
[S_{j}, t_\beta\S]_{*} \otimes S_{j}
 \rightarrow  t_\beta\S \rightarrow 
t_\alpha\S 
\end{equation}

Suppose that $\ell(w_{k})=1$, so that $w_{k}=s_{j}$ and $u=1$.  By induction on the number of Garside factors, we know that 
$[S_i,t_\beta \S]_{k+2}=0$, and $[S_j,t_\beta\S]_{k+2}=0$, 
and that $[S_i,[S_{j}, t_\beta\S]_{*} \otimes S_{j}]_{k+3}=0$, 
so (\ref{triangle2}) implies 
$[S_i,t_\alpha \S]_{k+2}=0$, as desired.  

\medskip

Now let $\ell(w_{k}) > 1$. We consider two cases.

\medskip

\noindent {\bf Case 1} Suppose $s_i$ and $s_j$ commute. Then
$u$ does not admit $s_i$ as a left factor, since we know that $s_ju$
does not admit $s_i$ as a left factor.  The induction hypothesis 
applied to $\beta$ implies that $[S_i,t_\beta \S]_{k+2}= 0$.  Apply
$t_j$ to both $S_i$ and $t_\beta\S$; since $t_jS_i \simeq S_i$, we
obtain that $[S_i,t_\alpha \S]_{k+2}=0$, as desired. 

\medskip

\noindent {\bf Case 2} Suppose $s_i$ and $s_j$ do not commute. 
Applying $\Hom(S_i,?)$ to (\ref{triangle2}) we get the following portion of a long exact sequence:

\begin{equation}\label{les} 
[S_{j}, t_\beta\S]_{k+1} \otimes [S_i,S_{j}]_{1}
 \rightarrow  [S_i,t_\beta\S]_{k+2} \rightarrow 
[S_i,t_\alpha\S]_{k+2}\rightarrow  
[S_{j}, t_\beta\S]_{k+2} \otimes [S_i,S_{j}]_{1}\end{equation}

Since $s_j$ is not a left factor of $u$, we see by the induction hypothesis applied to $t_\beta \S$ that the final term is zero.  
This implies that the map $ [S_i,t_\beta\S]_{k+2} \rightarrow [S_i,t_\alpha\S]_{k+2}$ is surjective, so 
any morphism in $[S_i, t_\alpha\S]_{k+2}$ factors through 
$t_\beta \S[k+2]$.  If $u$ does not admit $s_i$ as a left factor,
then we are done, since, by the induction hypothesis, $[S_i,t_\beta\S]_{k+2}=0$.
So suppose otherwise, and write $w_k=s_js_iv$.  Note that $v$ 
does not admit $s_j$ as a left factor, since $s_js_is_j=s_is_js_i$, 
which would imply that $w_k$ admits $s_i$ as a left factor.  
Let 
$\gamma=\tilde v\tilde{w}_{k-1} \cdots \tilde{w}_1$.  We now have:
$$[S_i,t_\gamma \S]_{*} \otimes S_i \rightarrow t_\gamma \S 
\rightarrow t_\beta \S$$

Since $v$ does not admit $s_{i}$ as a left factor, by the induction hypothesis
we know that $[S_i,t_\gamma\S]_{k+2}=0$, and thus any non-zero map
$\varphi: S_i \rightarrow t_\beta\S[k+2]$ composes to give a non-zero map 
$S_{i} \rightarrow [S_i,t_\gamma \S]_{*} \otimes S_i[k+3]$. Since $[S_{i},t_{\gamma}\S]_{d} = 0$ for $d > k+1$
by induction, such a map must land in the summand $[S_{i}, t_{\gamma}\S]_{k+1}\otimes S_{i}[2]$ of $[S_i,t_\gamma \S]_{*} \otimes S_i[k+3]$. By Lemma \ref{factor},
$\varphi$ must therefore factor through $S_{j}[1]$, so we have the arrows
in the diagram below other than the dotted arrow.  

$$\xymatrix{t_\beta\S[k+2] \ar[r]&  [S_i,t_\gamma]_{*} \otimes S_i[k+3] \ar[r] & 
t_{\gamma}\S[k+3]\\
             S_i \ar[u]\ar[r]&S_j[1]\ar[u]\ar@{.>}[ul] }$$

By the induction hypothesis (on length if $v \neq 1$ and on the number of Garside factors if $v=1$), we 
know that $\Hom(S_j[1],t_\gamma\S[k+3])=0$.  It follows that the
dotted arrow can be filled in, making the upper triangle commutative.  

We claim that the lower triangle must also be commutative. If it were not, the difference between
the two maps from $S_i$ to $t_\beta\S[k+2]$ would induce a non-zero map
to $t_\gamma\S[k+2]$, but no such map exists. 

Now consider (\ref{les}) again.  We have seen that any
$\varphi: S_i \rightarrow t_\beta\S[k+2]$ factors through $S_j[1]$,
so the map 
$[S_{j},t_{\beta}\S]_{k+1}\otimes [S_{i},S_{j}]_{1} \rightarrow [S_{i},t_{\beta}\S]_{k+2}$.
must be surjective. But we have already
argued that the map from $[S_i,t_\beta\S]_{k+2}$ to $[S_i,
t_\alpha\S]_{k+2}$ is surjective, so by exactness of (\ref{les}) we have $[S_i,t_\alpha\S]_{k+2}=0$, 
as desired.  
\end{proof}

\section{{\bf Application to spaces of stability conditions}}

We briefly recall the notion of stability condition on a triangulated category as introduced by Bridgeland in \cite{bridgeland}, review the results of Bridgeland in \cite{klein}, and point out how our faithfulness result Theorem \ref{faithful} answers a question left open in \cite{klein}.\\

A stability condition on a triangulated category
$\D$ consists of a bounded $t$-structure on $\D$ together with a  group homomorphism $\mathcal{Z}: K_{0}(\D) \rightarrow \mathbb{C}$, known as the  `central charge' such that
each object in the heart of the $t$-structure has a Harder-Narasimhan filtration with respect to the slope $Im(\mathcal{Z})/Re(\mathcal{Z})$ of the central charge. 
The set of stability conditions satisfying some technical hypotheses can be given the structure of a complex manifold ${\rm Stab}(\D)$.

Now consider the minimal resolution $\pi: X \rightarrow \mathbb{C}^{2}/G$ of the Kleinian 
singularity, as described in Example \ref{kleinian}. Inside the bounded derived category of coherent sheaves $D^{b}(X)$, let $\D$ be the thick subcategory generated by the $\Gamma$-configuration of $2$-spherical objects $S_{i}=\mathcal{O}_{E_{i}}(-1)[1]$. Bridgeland \cite{klein} has shown that there is a connected component 
${\rm Stab}_{0}(\D)$ of  ${\rm Stab}(\D)$, stable under the action of the braid group induced by the spherical twists $t_{i}$, together with a covering map $p:{\rm Stab}_{0}(\D) \rightarrow {\frak h}^{reg}/W$. As noted in Remark
\ref{fundamental}, it is well-known that the fundamental group $\pi_{1}({\frak h}^{reg}/W)$ 
is the braid group $B$ of the corresponding type. Bridgeland further shows that the image of $\rho: B \rightarrow {\rm Aut}(\D)$ gives the full group of deck transformations for the covering $p$, and therefore that if $\rho$ is injective, $p$
is in fact a universal cover. When Bridgeland 
was writing \cite{klein}, the injectivity of $\rho$ was known in type $A$ due to Seidel-Thomas \cite{st}. From Theorem \ref{faithful}, such injectivity is now known, and thus $p$ is a universal cover, in all types ADE.

\begin{remark}
In \cite{klein}, Bridgeland also considers an affine analogue of the category $\D$, with an extra generator $S_{0}= \mathcal{O}_{E}$, the structure sheaf of the exceptional divisor $\pi^{-1}(0)=E$, and establishes analogous results
for its space of stability conditions. 

It is therefore an interesting problem to determine if
the affine braid group action on this category (with extra generator given by the spherical twist along $S_{0}$) is faithful. This has been established by Ishii-Ueda-Uehara \cite{iuu} in type A. Somewhat more generally, one can ask
if the extended affine braid group action generated by
the finite type braid group and the Picard group of $X$ is
faithful.
\end{remark}

\bibliographystyle{plain}
\bibliography{bravthomasbiblio}

\medskip

\noindent Christopher Brav, University of Toronto, brav@math.toronto.edu

\medskip

\noindent Hugh Thomas, University of New Brunswick, hthomas@unb.ca

\end{document}